\documentclass[10pt,twoside]{article}
\usepackage{graphicx}
\usepackage{amsmath}
\usepackage{Latex-document}

\title{\bf  Some Aspects of Additive Coalescents\vskip 6mm}
\author{Jean Bertoin\vspace*{-0.5cm}\thanks
{Laboratoire de Probabilit\'es et Mod\`eles Al\'eatoires  and
Institut universitaire de France, Universit\'e Paris 6, 175, rue
du Chevaleret,  F-75013 Paris, France. E-mail:
jbe@ccr.jussieu.fr}}
\date{\vspace{-8mm}}

\begin{document}

\maketitle

\thispagestyle{first} \setcounter{page}{15}

\begin{abstract}

\vskip 3mm

We present some aspects of the so-called additive coalescence,
with a focus on its connections with random trees, Brownian
excursion, certain  bridges with exchangeable increments, L\'evy
processes, and sticky particle systems.

\vskip 4.5mm


\end{abstract}

\vskip 12mm

\section{ Additive coalescence in the finite setting} \label{section
1}\setzero \vskip-5mm \hspace{5mm}

The additive coalescence is a simple Markovian model for random
aggregation that arises for instance in the study of droplet
formation in clouds \cite{Go, D}, gravitational clustering in the
universe \cite{SP}, phase transition for parking \cite{CL},  ...
As long as only finitely many clusters are involved, it can be
described as follows. A typical configuration is a finite sequence
$x_1\geq \ldots\geq  x_n>0$ with $\sum_{1}^{n}x_i=1$, which may be
thought of as the ranked sequence of masses of clusters in a
universe with unit total mass. Each pair of clusters, say with
masses $x$ and $y$, merges as a single cluster with mass $x+y$ at
rate $K(x,y)=x+y$, independently of the other pairs in the system.
This means that to each pair $(i,j)$ of indices with $1\leq i< j
\leq n$, we associate an exponential variable $e(i,j)$ with
parameter $x_i+x_j$, such that to different pairs correspond
independent variables. If the minimum $\gamma_1:=\min_{1\leq i < j
\leq n}e(i,j)$ of these variables is reached, say for the pair
$(i_0,j_0)$, i.e. $\gamma_1=e(i_0,j_0)$, then at time $\gamma_1$,
we replace the clusters with labels $i_0$ and $j_0$  by a single
cluster with mass $x_{i_0}+x_{j_0}$. Then the system keeps
evolving with the same dynamics until it is reduced to a single
cluster.

An additive coalescent $(X(t), t\geq0)$ started from a finite
number $n$ of masses is a Markov chain in continuous times for
which the sequence of jump times $\gamma_1<\ldots<\gamma_{n-1}$
has a simple structure. Specifically, the increments between
consecutive coalescence times,
$\gamma_1,\gamma_2-\gamma_1,\ldots,\gamma_{n-1}-\gamma_{n-2}$, are
independent exponential variables with parameters $n-1,n-2,
\ldots, 1$, and are independent of the state chain
$\left(X(\gamma_k), k=0, \ldots,n-1\right)$. This elementary
observation enables one to focus on the chain of states, and to
derive simple relations with other random processes, as we shall
now see.

Pitman \cite{Pi} pointed at the following connection with random
trees. Pick a tree $\tau^{(n)}$ at random, uniformly amongst the
$n^{n-2}$ trees on $n$ labelled vertices.  Enumerate its $n-1$
edges at random and decide that they are all closed at the initial
time. At time $k=1,\ldots, n-1$, open the $k$-th edge, and  say
that two vertices belong to the same sub-tree if all the edges of
the path connecting those two vertices are open at time $k$.  If
we denote by $nY^{(n)}(k)$ the ranked sizes of sub-trees at time
$k$, then the chain $Y^{(n)}=(Y^{(n)}(k), k=0,\ldots n-1)$ has the
same law as the chain of states $(X^{(n)}(\gamma_k),k=0,\ldots
n-1)$ of the additive coalescent started from the initial
configuration $(1/n,\ldots,1/n)$.

\begin{picture}(300,60)(10,20)
\multiput(30,20)(50,0){7}{$\diamond$}
\multiput(40,20)(10,0){4}{$\cdot$}
\multiput(140,20)(10,0){4}{$\cdot$}
\multiput(240,20)(10,0){4}{$\cdot$}

\put (85,22){\line(1,0){45}} \put (185,22){\line(1,0){45}} \put
(285,22){\line(1,0){45}} \put (282,25){\line(0,1){35}} \put
(82,25){\line(0,1){35}} \put (82,25){\line(-2,3){25}} \put
(82,25){\line(2,3){25}} \put (133,25){\line(0,1){35}} \put
(282,25){\line(-2,3){25}} \put (282,25){\line(2,3){25}}

\put (55,62){\makebox(0,0){$\diamond$}}

\put (83,62){\makebox(0,0){$\diamond$}} \put
(108,62){\makebox(0,0){$\diamond$}} \put
(133,62){\makebox(0,0){$\diamond$}} \put
(255,62){\makebox(0,0){$\diamond$}} \put
(283,62){\makebox(0,0){$\diamond$}} \put
(308,62){\makebox(0,0){$\diamond$}}

\end{picture}
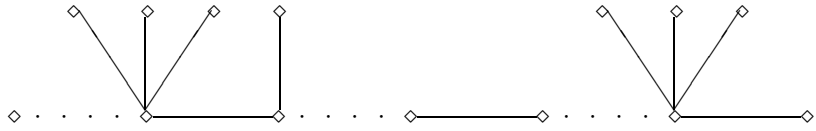
\vskip 10pt \centerline {\sl Forest derived at time 10 from a
tree with 14 vertices}

\vskip 10pt

Second, we lift from \cite{Be3} a different construction which is
closely related to hashing with linear probing, cf.
 \cite{CL}. We view the initial configuration
$x_1\geq \ldots\geq x_n>0$ as the ranked jumps of some bridge
$b=(b(u), 0\leq u \leq 1)$ with exchangeable increments. That is
we introduce $U_1,\ldots,U_n$, $n$ independent and uniformly
distributed variables and define
\begin{equation}\label{eq1}
b(u)\,=\,\sum_{i=1}^{n}x_i\left({\bf 1}_{\{u\geq
U_i\}}-u\right)\,,\qquad 0\leq u \leq 1\,.
\end{equation}

\begin{picture}(200,80)(-80,20)
\put (0,40){\line(1,0){180}} \put (0,40){\line(1,-1){50}}
\put(50,0){\vector(0,-1){10}} \put (50,10){\makebox(0,0){$r_3$}}
\put(50,20){\vector(0,1){10}} \put (50,30){\line(1,-1){50}}
\put(100,5){\vector(0,-1){25}} \put (100,18){\makebox(0,0){$r_1$}}
\put(100,25){\vector(0,1){45}} \put (100,70){\line(1,-1){50}}
\put(150,45){\vector(0,-1){25}} \put
(150,50){\makebox(0,0){$r_2$}} \put(150,55){\vector(0,1){15}}
\put (150,70){\line(1,-1){30}} \put (90,50){\makebox(0,0){$\mu$}}
\put (0,50){\makebox(0,0){$0$}} \put (180,50){\makebox(0,0){$1$}}
\put (100,40){\makebox(0,0){$\times$}} \put
(0,40){\makebox(0,0){$\times$}} \put
(180,40){\makebox(0,0){$\times$}} \put
(50,40){\makebox(0,0){$\times$}} \put (50,50){\makebox(0,0){${
}_{U_3}$}}
\end{picture}

\vskip50pt \centerline{\sl bridge $b$} \vskip 5mm \noindent

Next, we consider a path transformation (see the picture below)
that has been introduced by Tak\' acs \cite{Ta} and used by
Vervaat \cite{ver} to change a Brownian bridge into a normalized
Brownian excursion. Specifically, we set
\begin{equation}\label{eq2}
\epsilon(u)\,=\,b(u+\mu\,{\rm [mod\  1]})-b(\mu-)\,,\qquad 0\leq
u \leq 1\,,
\end{equation}
where $\mu$ stands for the location of the infimum of the bridge
$b$.

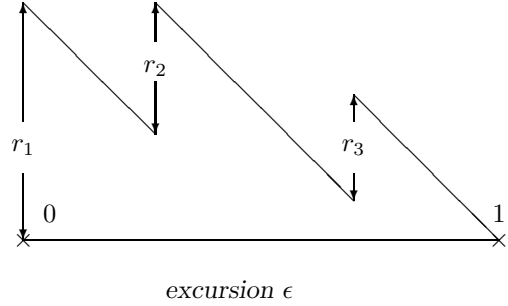
\begin{picture}(200,80)(120,20)
\put (200,-20){\line(1,0){180}} \put(200,5){\vector(0,-1){25}}
\put (200,15){\makebox(0,0){$r_1$}} \put(200,25){\vector(0,1){45}}
\put (200,70){\line(1,-1){50}} \put(250,40){\vector(0,-1){20}}
\put (250,45){\makebox(0,0){$r_2$}} \put(250,55){\vector(0,1){15}}
\put (250,70){\line(1,-1){75}} \put (325,35){\line(1,-1){55}}
\put(325,5){\vector(0,-1){10}} \put (325,15){\makebox(0,0){$r_3$}}
\put(325,25){\vector(0,1){10}}

\put (210,-10){\makebox(0,0){$0$}} \put
(380,-10){\makebox(0,0){$1$}} \put
(200,-20){\makebox(0,0){$\times$}} \put
(380,-20){\makebox(0,0){$\times$}}
\end{picture}

\vskip50pt \centerline{\sl excursion $\epsilon$} \vskip 5mm
\noindent

Finally, for every $t\geq0$, call $t$-interval any maximal
interval $[a,b[\subseteq [0,1]$ on which
$$tu-\epsilon(u)\,<\,\max\left\{(tv-\epsilon(v))^+,  0\leq v\leq
u\right\}\,, \qquad \hbox{for all } u\in[a,b]\,.
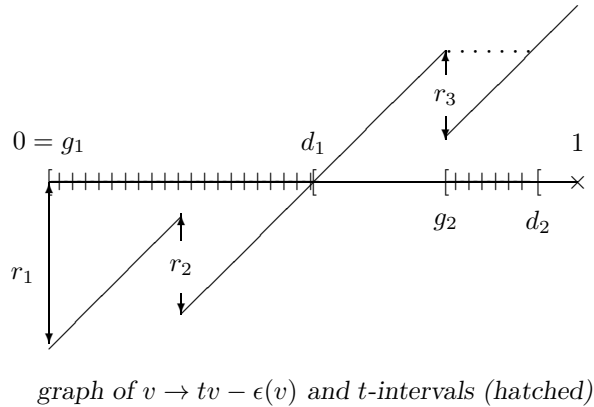
$$
\begin{picture}(300,110)(-30,20)
\put (50,50){\line(1,0){200}} \put (50,-13){\line(1,1){50}}
\put(50,15){\vector(0,1){35}} \put (40,15){\makebox(0,0){$r_1$}}
\put(50,24) {\vector(0,-1){35}}

\put (100,0){\line(1,1){100}} \put(100,8){\vector(0,-1){8}} \put
(100,18){\makebox(0,0){$r_2$}} \put(100,28) {\vector(0,1){10}}

\put (200,67){\line(1,1){50}} \put(200,91){\vector(0,1){9}} \put
(200,82){\makebox(0,0){$r_3$}} \put(200,75) {\vector(0,-1){9}}

\put (50,50){\makebox(0,0){$[$}} \put
(250,50){\makebox(0,0){$\times$}} \put
(50,65){\makebox(0,0){$0=g_1$}} \put
(150,65){\makebox(0,0){$d_1$}} \put (250,65){\makebox(0,0){$1$}}
\put (200,35){\makebox(0,0){$g_2$}} \put
(150,50){\makebox(0,0){$[$}} \put (200,50){\makebox(0,0){$[$}}
\multiput(200,48)(5,0){6}{$+$} \multiput(200,97)(5,0){7}{$\cdot$}
\multiput(50,48)(5,0){20}{$+$} \put (235,35){\makebox(0,0){$d_2$}}
\put (235,50){\makebox(0,0){$[$}}
\end{picture}

\vskip40pt \centerline {\sl graph of $v\to tv-\epsilon(v)$ and
$t$-intervals (hatched)} \vskip 5mm It is easy to see that the
$t$-intervals get finer as $t$ increases and tend to reduce to
the jump times of $\epsilon$ when $t\to\infty$. Denote by $F(t)$
the ranked sequence of the sums of the jumps made by $\epsilon$
on each $t$-interval, and by $0<\delta_1<\ldots <\delta_{n-1}$
the jump times of $F(\cdot)$. Then the chain $(F(\delta_{n-k-1}),
k=0,\ldots n-1)$ has the same law as the chain of states
$(X(\gamma_k),k=0,\ldots n-1)$ of the additive coalescent started
from the initial configuration $(x_1,\ldots,x_n)$.

\section{Standard and other eternal coalescents} \label{section 2}
\setzero\vskip-5mm \hspace{5mm}

Dealing with a finite number of clusters may be useful to give a
simple description of the dynamics, however it is a rather
inconvenient restriction in practice. In fact, it is much more
natural to work with  the infinite simplex
$${\cal S}^{\downarrow}\,=\,\left\{x=(x_1,x_2,\ldots): x_i\geq0 \hbox{
and }\sum_{i=1}^{\infty}x_i=1\right\}$$ endowed with the uniform
distance. In this direction, Evans and Pitman \cite{EP} have
shown that the semigroup of the additive coalescence enjoys the
Feller property on ${\cal S}^{\downarrow}$. Approximating a
general configuration $x\in{\cal S}^{\downarrow}$ by
configurations with a finite number of clusters then enables us
to view the additive coalescence as a Markovian evolution on
${\cal S}^{\downarrow}$. It is interesting in this setting to
consider asymptotics when the coalescent starts with a large
number of small clusters, which we shall now discuss.

Evans and Pitman \cite{EP} have first observed that the so-called
{\it standard additive coalescent} $(X^{(\infty)}(t),-\infty<
t<\infty)$ arises at the limit as $n\rightarrow\infty$ of the
additive coalescent process $(X^{(n)}(t),  -\frac{1}{2} \log n\leq
t<\infty)$ started at time $-\frac{1}{2} \log n$ with $n$
clusters, each with mass $1/n$. This limit theorem is perhaps
better understood if we recall the connection with the uniform
random tree $\tau^{(n)}$ on $n$ vertices that was presented in the
previous section. Indeed, if one puts a mass $1/n$ at each vertex
and let each edge have length $n^{-1/2}$, then $\tau^{(n)}$
converges weakly as $n\to\infty$ to the so-called {\it continuum
random tree} $\tau^{(\infty)}$; see Aldous \cite{A0}. More
precisely, $\tau^{(\infty)}$ is a compact metric space endowed
with a probability measure (arising as the limit of the masses on
vertices) which is concentrated on the leaves of the tree, and a
skeleton equipped with a length measure which is used to define
the distance between leaves. This suggests that the standard
additive coalescent might be constructed as follows: as time
passes, one creates a continuum random forest by logging the
continuum random tree along its skeleton and consider the ranked
sequence of masses of the subtrees. This yields a fragmentation
process, and the standard additive coalescent is finally obtained
by time-reversing this fragmentation process. Aldous and Pitman
\cite{AP1} have made this construction rigorous; more precisely
they showed that the tree $\tau^{(\infty)}$ has to be cut at
points
 that appear according to a Poisson point process on the
skeleton with intensity given by the length measure. This
representation yields a number of explicit  statistics for the
standard additive coalescent. For instance, for every $t\in{\bf
R}$, the distribution of $X^{(\infty)}_t$ is given by that of the
ranked sequence $\xi_1\geq\xi_2\geq\ldots$ of the atoms of a
Poisson measure on $]0,\infty[$ with intensity ${\rm e}^{-t}(2\pi
x^3)^{-1/2}dx$ and conditioned by $\xi_1+\cdots=1$.

The continuum random tree bears remarkable connections with the
{\it Brownian excursion} (cf. for instance Le Gall \cite{LG}),
and one naturally expects that the standard additive coalescent
could also be constructed from the latter. This is indeed feasible
(see \cite{Be2} and also \cite{CL}) although its does not seem
obvious to relate the following construction with that based on
the continuum random tree. Specifically, let $(\epsilon(s), 0\leq
s \leq 1)$ be a Brownian excursion with unit duration, and for
every $t\geq0$, consider the random open set
\begin{equation}\label{eq3}
G(t)\,=\,\left\{s\in[0,1]: ts-\epsilon(s) < \max_{0\leq u \leq
s}\left(tu-\epsilon(u)\right)\right\}\,.
\end{equation}
Then $G(t)$ decreases as $t$ increases, and if we denote by $F(t)$
the ranked sequence of the lengths of its intervals components
(which of course are related to the so-called $t$-intervals of the
preceding section), then $(F({\rm e}^{-t}), -\infty<t<\infty)$ is
a standard additive coalescent.

More generally, Aldous and Pitman \cite{AP2} have characterized
all the processes that may arise as the limit of additive
coalescents started with a large number of small clusters. They
are referred to as {\it eternal} additive coalescents as these
processes are indexed by times in $]-\infty,\infty[$. They can be
constructed by the same procedure as in \cite{AP1} after
replacing the continuum random tree $\tau^{(\infty)}$ by a
so-called inhomogeneous continuum random tree.

An alternative construction was proposed in \cite{Be3} and
\cite{Mie}. Specifically, one may replace the standard Brownian
excursion $\epsilon$ by that obtained by the Tak\'acs-Vervaat
transformation (\ref{eq2}) where $(b(s), 0\leq s \leq 1)$ is now a
bridge with exchangeable increments, no positive jumps and
infinite variation (which arises as the limit of elementary
bridges of the type (\ref{eq1}), see Kallenberg \cite{Kal}).  The
ranked sequence $F(t)$ of the lengths of the interval components
of $G(t)$ defined by (\ref{eq3}) then yields a fragmentation
process, and by time-reversal, $(F({\rm e}^{-t}),
-\infty<t<\infty)$ is an eternal additive coalescent.

Roughly, this construction can be viewed as the limit of that
presented in Section 1 when the additive coalescent starts from a
finite number of clusters.

\section{Eternal coagulation and certain L\'evy processes}
\label{section 3} \setzero\vskip-5mm \hspace{5mm}

A long time before the notion of stochastic coalescence was
introduced, Smoluchowski \cite{Sm} considered a family of
differential equations to model the evolution in the hydrodynamic
limit of a particle system in which particles coagulate pairwise
as time passes.  It bears natural connections with the stochastic
coalescence; we refer to the survey by Aldous \cite{A2} for
detailed explanations, physical motivations, references ...
Typically, we are given a symmetric kernel $K:]0,\infty[\times
]0,\infty[\to[0,\infty[$ that specifies the rate at which two
particles coagulate as a function of their masses. Here, we take
of course $K(x,y)=x+y$. If we represent the density of particles
with mass ${\tt d}x$ at time $t$ by a measure $\mu_t({\tt d}x)$
on $]0,\infty[$, then
\begin{equation}\label{eqSm}
\frac{\tt d}{{\tt d}t}\langle\mu_t,f\rangle \,=\, \frac 1 2
\int_{]0,\infty[\times]0,\infty[}
\left(f(x+y)-f(x)-f(y)\right)(x+y)\mu_t({\tt d}x)\mu_t({\tt
d}y)\,,
\end{equation}
where $f$ a test function and $\langle\mu_t,f\rangle=\int
f(x)\mu_t({\tt d}x)$. Motivated by the preceding section, we are
interested in the eternal solutions of (\ref{eqSm}), in the sense
that the time parameter $t$ is real (possibly negative). It is
proven in \cite{Be4} that every eternal solution
$(\mu_t)_{t\in{\bf R}}$ subject to the normalizing condition
$\int_{}^{}x\mu_t({\tt d}x)=1$ (i.e. the total mass of the system
is $1$), can be constructed as follows.

First, define the function
\begin{equation}\label{eq5}
\Psi_{\sigma^2,\Lambda}(q)\,=\, \frac 1 2 \sigma^2
q^2+\int_{]0,\infty[}\left({\rm e}^{-qx}-1+qx\right)\Lambda({\tt
d}x)\,,\qquad q\geq0\,,
\end{equation}
where $\sigma^2>0$ and  $\Lambda$ is a measure on $]0,\infty[$
with $\int (x\wedge x^2)\Lambda({\tt d}x)<\infty$. We further
impose that either $\sigma^2>0$ or $\int x\Lambda({\tt
d}x)=\infty$. Next, let $\Phi(\cdot,s)$ be the inverse  of the
bijection $q\to\Psi_{\sigma^2,\Lambda}(sq)+q$. One can check that
$\Phi(q,{\rm e}^t)$ can be expressed in the form
\begin{equation}\label{eq6}
\Phi(q,{\rm e}^t)\,=\,\int_{]0,\infty[}(1-{\rm e}^{-qx})\mu_t({\tt
d}x)\,,\qquad q\geq0\,,
\end{equation} where $(\mu_t)_{t\in{\bf R}}$ is then an eternal solution
to Smoluchowski's coagulation equation. For instance, when
$\sigma^2=1$ and $\Lambda=0$, $\xi$ is a standard Brownian motion
and we recover the well-known solution
$$
\mu_t({\tt d}x)\,=\,\frac{{\rm e}^{-t}}{ \sqrt{2\pi
x^3}}\,\exp\left(-\frac{x{\rm e}^{-2t}}{2}\right){\tt d}x\,,\qquad
t\in{\bf R},\, x>0\,.$$

This invites a probabilistic interpretation. Indeed, (\ref{eq5})
is a special kind of {\it L\'evy-Khintchine formula}; see section
VII.1 in \cite{Be}. More precisely, there exists a {\it L\'evy
process with no positive jumps}, $\xi=\left(\xi_r, r\geq
0\right)$, such that
$${\bf E}\left(\exp\left(q\xi_r\right)\right)\,=\,\exp
\left(r\Psi_{\sigma^2,\Lambda}(q)\right)\,,\qquad q\geq0\,.$$ It
is then well-known (e.g.  Theorem VII.1 in \cite{Be}) that the
first passage process
$$T^{(s)}_x\,:=\,\inf\left\{r\geq0: s\xi_r+r>x\right\}\,,\qquad x\geq0$$
is a subordinator with
$${\bf E}\left(\exp\left(-qT^{(s)}_x\right)\right)\,=\,\exp
\left(-x\Phi(q,s)\right)\,,\qquad q,x\geq0\,,$$ where the Laplace
exponent $\Phi(\cdot,s)$ is the inverse bijection of
$q\to\Psi_{\sigma^2,\Lambda}(sq)+q$. Thus (\ref{eq6}) can be
interpreted as the L\'evy-Khintchine formula for $\Phi(\cdot,s)$,
and we conclude that the eternal solution $\mu_t$  can be
identified as the  L\'evy measure of the subordinator $T^{(s)}$
for $s={\rm e}^t$.

This probabilistic interpretation also points at a simple random
model for aggregation of intervals. Indeed, the closed range
${\cal T}^{(s)}=\left\{T^{(s)}_x, x\geq0\right\}^{\rm cl}$ of
$T^{(s)}$ induces a partition of $[0,\infty[$ into a family of
random disjoint open intervals, namely the interval components of
$G(s)=[0,\infty[\backslash {\cal T}^{(s)}$. We now make the key
observation that
\begin{equation}\label{eqz}
{\cal T}^{(s)}\subseteq {\cal T}^{(s')}\qquad \hbox{for }0<s'<s,
\end{equation}
because an instant at which $r\to s\xi_r+r$ reaches a new maximum
is always also an instant at which $r\to s'\xi_r+r$ reaches a new
maximum. Roughly, (\ref{eqz}) means that the random partitions
get coarser as the parameter $s$ increases; and therefore they
induce a process in which intervals aggregate. The latter is
closely related to a special class of eternal additive
coalescents, and has been studied in \cite{Be2, Sch} in the
Brownian case, and in \cite{Mie} in the general case.

\section{Sticky particle systems} \label{section 4}
\setzero\vskip-5mm \hspace{5mm}

Sticky particle systems evolve according to the dynamics of
completely inelastic collisions with conservation of mass and
momentum, which are also known as the dynamics of ballistic
aggregation. This means that the velocity of particles only
changes in case of collision, and in that case, a heavier cluster
merges at the location of the shock with mass and momentum given
by the sum of the masses and momenta of the clusters involved.
This has been proposed as a model for the formation of large scale
structures in the universe; see the survey article \cite{V}. We
now have two quite different dynamics for clustering: on the one
hand the ballistic aggregation which is deterministic, and on the
other hand the additive coalescence which is random and may appear
much more elementary as it does not take into account significant
physical parameters such as distances between clusters and the
relative velocities. Nonetheless, there is a striking connection
between the two when randomness is introduced in the deterministic
model, as we shall now see.

We henceforth focus on dimension one and assume that at the
initial time, particles are infinitesimal (i.e. fluid) and
uniformly distributed on the line. The evolution of the sticky
particle system can then be completely analyzed in terms of the
entropy solution to a single PDE, the transport equation
\begin{equation}\label{eqburgers}
\partial_t u+u\partial_x u\,=\,0\,.
\end{equation}
Here $u(x,t)$ represents the velocity of the particle located at
$x$ at time $t$, and the entropy condition imposes that for every
fixed $t>0$, the function $u(\cdot,t)$ has only discontinuities of
the first kind and no positive jumps (the latter restriction
accounts for the total inelasticity of collisions). Provided that
the initial velocity $u(\cdot,0)$ satisfies some very mild
hypothesis on its rate of growth,  there is a unique weak
solution to the equation (\ref{eqburgers}) which fulfills the
entropy condition, and which can be given explicitly in terms of
$u(\cdot,0)$.

 We assume that the initial velocities in the particle system are
random, and more precisely
$$u(r,0)=0 \hbox{ for }r<0 \quad \hbox{and} \quad (u(r,0), r\geq0)
 \mbox{$ \ \stackrel{\cal L}{=}$ }(\xi_r, r\geq0)\,,$$
where $\xi$ denotes the L\'evy process with no positive jumps
which was used in the preceding section.

Roughly, the  dynamics of sticky particles are not only
deterministic, but also induce a loss of information as time goes
by, in the sense that the initial state of the system entirely
determines the state at time $t>0$, but cannot be completely
recovered from the latter. In this direction, let us observe the
system at some fixed time $t>0$, i.e. we know the locations,
masses and velocities of the clusters at this time. Let us pick a
cluster located in $[0,\infty[$, using for this only the
information available at time $t$ (for instance, we may choose
the heaviest cluster located at time $t$ in $[0,1]$). We shall
work conditionally on the mass of this cluster, and for
simplicity, let us assume it has unit mass. For every
$r\in]0,t[$, denote by $M(r)=\left( m_1(r),  m_2(r),
\ldots\right)$ the ranked sequence of masses of clusters at time
$r$ which, by time $t$ have aggregated to form the cluster we
picked, so $M(r)$ can be viewed as a random variable with values
in ${\cal S}^{\downarrow}$. Then the time-changed processes
$$M\left(t\left(1- \frac{t}{t+{\rm e}^s}\right)\right)\,,\qquad
-\infty<s<\infty $$
 is an eternal additive coalescent.
This was established in \cite{Be1} in the case of Brownian initial
velocity; and the recent developments on eternal additive
coalescents made in \cite{AP2, Be3, Mie} show that the arguments
also applies for L\'evy type initial velocities.

\label{lastpage}

\end{document}